\documentclass[leqno]{article}
 
\usepackage[papersize={176.76mm,250mm},margin=0.75in ]{geometry}

\usepackage{etex}
\usepackage{amsmath, bm,amssymb,amsfonts,dsfont,amsthm}
\usepackage{graphicx}
\usepackage[usenames,dvipsnames,svgnames,table]{xcolor}
\usepackage{multicol}
\usepackage{mathrsfs}
\usepackage{mathtools}
\usepackage{amsthm}
\usepackage{enumitem}
\usepackage{lineno}
\definecolor{antiquefuchsia}{rgb}{0.57, 0.36, 0.51}
\definecolor{auburn}{rgb}{0.43, 0.21, 0.1}
\definecolor{darkcerulean}{rgb}{0.03, 0.27, 0.49}
\definecolor{denim}{rgb}{0.08, 0.38, 0.74}
\definecolor{black}{rgb}{0.0, 0.0, 0.0}
\definecolor{sacramentostategreen}{rgb}{0.0, 0.34, 0.25}
\definecolor{phthaloblue}{rgb}{0.0, 0.06, 0.54}
\usepackage[colorlinks=true,linkcolor=black, urlcolor=auburn, filecolor=blue, citecolor=sacramentostategreen,backref=page]{hyperref}

\usepackage{etoolbox}
\makeatletter
\patchcmd{\BR@backref}{\newblock}{\newblock(Cited on page~}{}{}
\patchcmd{\BR@backref}{\par}{)\par}{}{}
\makeatother

\usepackage{relsize}

\usepackage[all]{xy}
\usepackage{tikz-cd}
\usepackage{array}
\usepackage{tensor}
\usepackage[cal=mathpi,calscaled=.94, bb=ams,frakscaled=.97,scr=rsfso]{mathalfa} 
\usepackage[utf8]{inputenc}

\usepackage{pdflscape}
\usepackage{float}

\usepackage{scalerel,stackengine}

\usetikzlibrary{matrix}
\usetikzlibrary{arrows}

\newcommand{\ds}[1]{\ensuremath{\mathds{#1}}}

\newcommand{\curly}[1]{\ensuremath{\mathscr{#1}}}

\newcommand{\prj}{\ensuremath{\mathrm{Proj}}}

\DeclareMathOperator{\bl}{Bl}

\newcommand{\perfp}{\ensuremath{\ds{P}_K^{n,\mathrm{ad,perf}}}}
\newcommand{\perfpa}[1]{\ensuremath{\ds{P}_K^{{#1},\mathrm{ad,perf}}}}

\newcommand{\pro}{\ensuremath{\ds{P}_K^{1,\mathrm{ad,perf}}}}

\newcommand{\aro}{\ensuremath{\ds{A}_K^{1,\mathrm{ad,perf}}}}

\newcommand{\arov}[1]{\ensuremath{\ds{A}_K^{{#1},\mathrm{ad,perf}}}}

\newcommand{\lr}[1]{\left\langle{#1}\right\rangle}

\newcommand{\id}[1]{\ensuremath{{\mathfrak{#1}}}}

\newcommand{\abs}[1]{\ensuremath{{\left\vert{#1}\right\vert}}}

\newcommand{\xra}[1]{\xrightarrow{#1}}
\newcommand{\ra}{\rightarrow}
\newcommand{\til}[1]{\widetilde{#1}}

\newcommand{\ten}[4]{\ensuremath{\tensor[^{#1}]{#2}{_{#3,#4}}}}
\usepackage{ccfonts,eucal}
\usepackage[euler-digits,euler-hat-accent]{eulervm}
\usepackage{sseq}

\newtheorem{theorem}{Theorem}[section]

\newtheorem{mdef}[theorem]{Definition}
\let\olddefinition\mdef
\renewcommand{\mdef}{\olddefinition\normalfont}

\newtheorem{exam}[theorem]{Example}
\let\oldexample\exam
\renewcommand{\exam}{\oldexample\normalfont}

\newtheorem{rem}{Remark}
\let\oldremark\rem
\renewcommand{\rem}{\oldremark\normalfont}

\stackMath
\newcommand\rwhat[1]{%
\savestack{\tmpbox}{\stretchto{%
  \scaleto{%
    \scalerel*[\widthof{\ensuremath{#1}}]{\kern-.6pt\bigwedge\kern-.6pt}%
    {\rule[-\textheight/2]{1ex}{\textheight}}
  }{\textheight}%
}{0.5ex}}%
\stackon[1pt]{#1}{\tmpbox}%
}
\usepackage{changepage}
\usepackage{fancyhdr} 
\newcommand{\myreferences}{C:/Users/Harpreet/Documents/Biblio/bibThesis}

\begin{document}

\title{Geometry of Projective Perfectoid and Integer Partitions}
\author{Harpreet Singh Bedi~~~{bedi@gwu.edu}}

\maketitle

\begin{abstract}
Line bundles of rational degree are defined using Perfectoid spaces, and their cohomology computed via standard \v{C}ech complex along with Kunneth formula. A new concept of `braided dimension' is introduced, which helps convert the curse of infinite dimensionality into a boon, which is then used to do Bezout type computations, define euler characters, describe ampleness and link integer partitions with geometry. This new concept of 'Braided dimension' gives a space within a space within a space an infinite tower of spaces, all intricately braided into each other.  Finally, the concept of Blow Up over perfectoid space is introduced.
\end{abstract}

\subsection*{Acknowledgement} I am grateful to Kiran Kedlaya for critical comments. Any remaining errors are my own. Many thanks to Yongwu Rong, Jozef Przytycki, Robbie Robinson, Murli Gupta and Niranjan Ramachandran for help and support. The work in this paper was part of \cite{Bedi2018}, and some concepts have been further built upon.

\tableofcontents

\section{Introduction}

We begin the story by recalling the notion of degree in topology as given in \cite[pp 134]{hatcher2002algebraic}. A map of spheres $f:S^n\ra S^n$ induces a map in homology $\phi:H_n(S^n)\ra H_n(S^n)$ given as a group homomorphism $\ds{Z}\ra \ds{Z}$. This map is simply multiplication by $d\in\ds{Z}$ which is called the degree of the map. 

If we take direct limit of these maps \eqref{eq:direct_lim} as in example 3F.3 \cite[pp 312]{hatcher2002algebraic} by setting $d=p$ a prime

\begin{equation}\label{eq:direct_lim}
\ds{Z}\xra{p}\ds{Z}\xra{p}\ds{Z}\ra\cdots
\end{equation}
we get a moore space $M(\ds{Z}[1/p],n)$. We can now talk about maps from $\ds{Z}[1/p]\ra \ds{Z}[1/p]$, and degree $d$ as an element of $\ds{Z}[1/p]$. In this article we transfer the notion of degree as an element of $\ds{Z}[1/p]$ to geometry by defining degree of line bundles on projective perfectoid space. The projective perfectoid space is obtained by gluing together affine perfectoid spaces and are constructed in \cite{scholze_1}. The foundational details can be found in \cite{Kiran_AWS} and \cite{Kiran_Liu_1}.

Recall that to a projective space we associate the ring of homogeneous polynomials $K[X_1,\ldots, X_n]$, the analogous situation for projective perfectoid space is to associate the ring of power series $K\lr{X_1^{1/p^\infty}, \ldots,X_n^{1/p^\infty} }$ which is the multivariate version of $K\lr{X^{1/p^\infty}}$ given as
\[K\lr{X^{1/p^\infty}}=\sum_{i\in \ds{Z}[1/p]_{\geq0}}a_iX^i,~~ a_i\in K, ~~\lim_{i\ra\infty}\abs{a_i}=0\]

At this point the reader should be thinking of above series as Tate Algebras with rational powers. The series $K\lr{X^{1/p^\infty}}$ lives on the perfectoid affine line \aro. We glue two affine lines together to get \perfpa{1} with associated laurent power series $K\lr{X^{1/p^\infty},X^{-1/p^\infty} }$. The same approach is used to construct \perfp.

\section{Ordering the elements}It is possible to write down the elements of $K\lr{X_1^{1/p^\infty}, \ldots,X_n^{1/p^\infty} }$ with a natural order. We describe one such order on $f\in K\lr{X^{1/p^\infty}} $. 

Recall that we can order rational numbers $a/b$ by considering them as a tuple $(a,b)$, we will use the same ordering. The terms of $f$ consist of $X^{a/p^b}$ with $a\in\ds{Z}_{>0}$ and $b\in\ds{Z}_{\geq0}$, we arrange the terms by ordering the tuple $(a,b)$ as in rational numbers. Notice that for $b=0$ we simply get an integer power $X^a$. 

\section{Proj of Graded Ring} 

Let $S=\oplus_{d}S_{d}$ be the graded homogeneous co-ordinate ring. Corresponding to the graded ring $S$ we have a sheaf of rings $\curly{O}_\prj(S)$ which will give the scheme $(\prj(S),\curly{O}_{\prj(S)})$, as given in ( \cite[\href{http://stacks.math.columbia.edu/tag/01M3}{Tag 01M3}]{stacks-project} or page 117 of \cite{hartshorne1977algebraic}).

Let $M=\oplus_{\ds{Z}}M_n$ be a graded $S$ module, that is $S_nM_m\subset M_{n+m}$ then there is a sheaf $\til{M}$ on the basis of standard open sets.
\begin{equation}
\til{M}(U)=M_{(f)}=\{mf^{-d}\in M_f\text{ such that }m\in M_{d\cdot\deg f}\}
\end{equation}

We want to construct a graded ring for fractional power series. 

\subsection{Grading} Analogous to the case of standard grading of homogeneous polynomials, we can construct a fractional grading with $d\in \ds{Z}[1/p]$ arranged in the canonical order. 
\begin{equation}
S_d:=\text{ Homogneous polynomials of degree } d
\end{equation}
\section{Defining $\curly{O}(m)$}
Let $B(n)$ denote the graded $S$-module defined by $B(n)_d=S_{n+d}$. This will be called twist of $B$. For a given $X=\prj (S)$,  let $\curly{O}_X(n)$ denote the $\curly{O}_X$ module $\til{B(n)}$. Let $f\in S$ be homogeneous of degree one (with affine open set $D_+(f)$ ), then we get

\begin{align}
B(n)_{(f)}&=f^nB_{(f)}\\
\curly{O}_X(n)\vert_{D_+(f)}&=f^n\curly{O}\vert_{D_+(f)}
\end{align}
Notice that $n\in\ds{Z}$, but we could choose $n\in\ds{Z}[1/p]$ which would give us rational degrees.
Informally, we can define $\curly{O}(m)$ with $m\in\ds{Z}[1/p]$ as 

\begin{equation}\label{eq:twist1}
\curly{O}(m)=\left\{\frac{f}{g}\text{ where }~f,g\in S \text{ and }\deg f-\deg g=m\right\}
\end{equation}

Although, the above construction seems natural we need to point out the approach in \cite{Das_2016} in particular Question 5.3.3 on Page 69:

``Does there exist a well defined notion of `degree' for line bundles on perfectoid space?''

We have answered the question here, but the reader needs to look at the above mentioned work to appreciate what has been made obvious here wasn't known before.

\subsubsection*{Remark}
\begin{enumerate}
\item Note that since, $f$ and $g$ are homogeneous, degree is well defined for the series (it is the same for all elements). 
\item $\Gamma(\pro\curly{O}(1))$ is generated by  $X,Y,\ldots,X^{q_i}Y^{q_j},\ldots$ where $q_i+q_j=1$. The basis reduces to $X,Y$ if we take exponent of $p$ in $X^a/p^b$ as $0$. Thus, $\dim\Gamma(\pro\curly{O}(1))=\infty$ and reduces to $2$ when we take exponent of $p$ to be $0$.
\item $\Gamma(\pro\curly{O}(1/p^i))$ is generated by  $X^{1/p^i},Y^{1/p^i},\ldots,X^{t_i}Y^{t_j},\ldots$ where $t_i+t_j=1/p^i$ and thus, $\dim\Gamma(\pro\curly{O}(1))=\infty$.
\end{enumerate}

\subsection{Twisting the sheaf $\curly{O}(m)$}
We have an isomorphism between graded modules, which is given as  

\begin{align}
\oplus_{d}S_d&\xra{\cdot S_n}\oplus_{d}S_{d}, ~~~~~~n\in\ds{Z}[1/p]\\
S_d&\mapsto S_{d+n}
\end{align}

We get the twist by tensoring with $\curly{O}(m)$ where $m\in \ds{Z}[1/p]$. The tensor products of sheaves $\curly{O}(m)\otimes\curly{O}(n)=\curly{O}(m+n)$ gives us an isomorphism with $\ds{Z}[1/p]$. 

\section{Cohomology of Line Bundles on \perfp}\label{coho1}

Recall that global sections $H^0(\ds{P}^n,\curly{O}_{\ds{P}^n}(m))$ are generated by homogeneous polynomials of degree $m$ in $n+1$ variables. For example, $H^0(\ds{P}^1,\curly{O}_{\ds{P}^1}(2))$ is generated by $x^2,xy,y^2$. We can similarly, define the global sections on $\pro$, of degree $2$ being generated by $x^{a_1/p^{n_1}}\cdot y^{a_2/p^{n_2}}$ such that 

\begin{equation}
\frac{a_1}{p^{n_1}}+\frac{a_2}{p^{n_2}}=2\text{ and }a_i\in\ds{Z}_{>0}, n_i\in \ds{Z}_{\geq0}
\end{equation}
Thus, we see that the dimension of $H^0(\pro,\curly{O}_{\pro}(2))$ is infinite, and $H^0(\ds{P}^1,\curly{O}_{\ds{P}^1}(2))\subset H^0(\pro,\curly{O}_{\pro}(2))$. For $n_i=0$ we recover the global sections  $H^0(\ds{P}^1,\curly{O}_{\ds{P}^1}(2))$.

\begin{align}
\dim  H^0(\perfp,\curly{O}_{\perfp}(m))&=\infty 
 \end{align}
 
The line bundles of the form $2/p$ can also be defined with infinitely many global sections $x^{2/p},x^{1/p}y^{1/p},y^{2/p}, x^{1/p^2}y^{2-1/p^2},\ldots $, showing that this is a richer world than the regular line bundles studied in algebraic geometry.

Recall, that we are working with non-archimedean spaces and thus using Grothendieck Topology, but luckily the cohomology can be computed using Cech complex and Leray's theorem for covers. Hence we use the arguments on page 225 of \cite{hartshorne1977algebraic} to get the first theorem. The proof below works for all $m\in\ds{Z}[1/p]$.
 
\begin{theorem}
\begin{enumerate}
\item $H^0(\perfp,\curly{O}_{\perfp}(m))$ is a free module of infinite rank.
\item $H^n(\perfp,\curly{O}_{\perfp}(-m))$ for $m>n $ is a free module of infinite rank.
\end{enumerate}
\end{theorem} 
 
 We take the standard cover of $\perfp$ by affine sets $\id{U}=\{U_i\}_{i}$ where each $U_i=D(x_i), i=0,\ldots,n$. 
 
 We get $H^0(\perfp,\mathcal {F})$ as the kernel of the following map
 \begin{equation}
 \prod S_{x_{i_0}}\ra \prod S_{x_{i_0}x_{i_1}}
  \end{equation}
  
  The element mapping to the Kernel has to lie in all the intersections $S=\cap_i S_{x_i}$, as given on \cite[pp 118]{hartshorne1977algebraic}.
  
 $H^n(\perfp,\mathcal{F})$ is the cokernel of the map 

\begin{equation}
d^{n-1}:\prod_k S_{x_0\cdots \hat{x}_k\cdots x_n}\ra S_{x_0\cdots x_n}
\end{equation}

$S_{x_0\cdots x_n}$ is a free $A$ module with basis $x_0^{l_0}\cdots x_n^{l_n}$ with each $l_i\in\ds{Z}[1/p]$. The image of $d^{n-1}$ is the free submodule generated by those basis elements with atleast one $l_i\geq 0$. Thus $H^n$ is the free module with basis as negative monomials

\begin{equation}
\{x_0^{l_0}\cdots x_n^{l_n}\}\text{ such that }l_i<0 
\end{equation}

The grading is given by $\sum l_i$ and there are infinitely many monomials with degree $-n-\epsilon$ where $\epsilon$ is something very small and $\epsilon\in\ds{Z}[1/p]$. Recall, that in the standard coherent cohomology there is only one such monomial $x_0^{-1}\cdots x_n^{-1}$ . For example, in case of $\ds{P}^2$ we have $x_0^{-1}x_1^{-1}x_2^{-1}$ but in $\perfpa{2}$ in addition to above we also have  $x_0^{-1/2}x_1^{-1/2}x_2^{-2}$. 

Recall that in coherent cohomology of $\ds{P}^n$ the dual basis of $x_0^{m_0}\cdots x_n^{m_n}$ is given by $x_0^{-m_0-1}\cdots x_n^{-m_n-1}$ and the operation of multiplication gives pairing. We do not have this pairing for $\perfp$, but we can pair $x_0^{m_0}$ with $x_0^{-m_0}$.

\begin{theorem}
$H^i(\perfp,\curly{O}_{\perfp})=0$ if $0<i<n$
\end{theorem}

We will use the proof from \cite[pp 474-475]{ravi2}, using the convention that $H^0$ denotes global sections. We will work with $\perfp, n=2$ for the sake of clarity, the case for general $n$ is identical. The \v{C}ech complex is given in \ref{check1}.

\begin{landscape}

\begin{figure}
\centering
 \begin{tikzpicture}
 []
        \matrix (m) [
            matrix of math nodes,
            row sep=0.5em,
            column sep=2.5em,
                   ]
{
 ~ & ~ &  |[name=qa1]|U_0 & |[name=qb1]|U_{0}\cap U_1 &~&~\\
 ~ & ~ & \oplus & \oplus &~&~\\
 |[name=qka]| 0 &|[name=qkb]| \ds{P}^2 & |[name=qkc]| U_1 &|[name=qkd]| U_{0}\cap U_2 & |[name=qke]| U_{0}\cap U_1\cap U_2 &|[name=qkf]|0\\
  ~ & ~ & \oplus & \oplus &~&~\\
 ~ & ~ & |[name=qa2]| U_2 & |[name=qb2]|U_1\cap U_2 &~&~\\
 ~ & ~ &  |[name=a1]|H^0\lr{X_0^{-p^{-\infty}}}& |[name=b1]|H^0\lr{X_0^{-p^{-\infty}},X_1^{-p^{-\infty}} }&~&~\\
 ~ & ~ & \oplus & \oplus &~&~\\
 |[name=ka]| 0 &|[name=kb]| H^0 & |[name=kc]|H^0\lr{X_1^{-p^{-\infty}}}&|[name=kd]|H^0\lr{X_0^{-p^{-\infty}},X_2^{-p^{-\infty}} }&|[name=ke]|H^0\lr{X_0^{-p^{-\infty}},X_1^{-p^{-\infty}},X_2^{-p^{-\infty}}}&|[name=kf]|0\\
  ~ & ~ & \oplus & \oplus &~&~\\
 ~ & ~ & |[name=a2]| H^0\lr{X_2^{-p^{-\infty}}}& |[name=b2]|H^0\lr{X_1^{-p^{-\infty}},X_2^{-p^{-\infty}} }&~&~\\
        };

        \path[overlay,->, font=\scriptsize,>=latex]
        (ka) edge (kb)
         (kb) edge (kc)
    
          (kd) edge (ke)
          (ke) edge (kf)
          
            (kb) edge (a1)
            (kb) edge (a2)
             (a1) edge (b1)
 (a1) edge (kd)
 (a2) edge (b2)
(a2) edge (kd)
(kc) edge (b1)
(kc) edge (b2)
(b1) edge (ke)
(b2) edge (ke)
 ;
 
  \path[overlay,<-, font=\scriptsize,>=latex]

        (qka) edge (qkb)
         (qkb) edge (qkc)
                    (qkd) edge (qke)
          (qke) edge (qkf)
          
            (qkb) edge (qa1)
            (qkb) edge (qa2)
             (qa1) edge (qb1)
 (qa1) edge (qkd)
 (qa2) edge (qb2)
(qa2) edge (qkd)
(qkc) edge (qb1)
(qkc) edge (qb2)
(qb1) edge (qke)
(qb2) edge (qke)
 ;
\end{tikzpicture}
\caption{\v{C}ech Complex for $\perfp,n=2$}\label{check1}
\end{figure}

\end{landscape}
\begin{description}
\item[$3$ negative exponents] The monomial $X_0^{a_0}\cdot X_1^{a_1}\cdot X_2^{a_2}$ where $a_i<0$. We cannot lift it to any of the coboundaries (that is lift only to $0$ coefficients). If $K_{012}$ denotes the coefficient of the monomial in the complex (Figure \ref{check2}), we get zero cohomology except for the spot corresponding to $U_0\cap U_1\cap U_2$.

\begin{figure}[H]
\centering
 \begin{tikzpicture}
 []
        \matrix (m) [
            matrix of math nodes,
            row sep=0.5em,
            column sep=2.5em,
                   ]
{ ~ & ~ &  |[name=a1]|0 & |[name=b1]|0 &~&~\\
 ~ & ~ & \oplus & \oplus &~&~\\
 |[name=ka]| 0 &|[name=kb]| 0 & |[name=kc]|0 &|[name=kd]|0 &|[name=ke]|K_{012}&|[name=kf]|0\\
  ~ & ~ & \oplus & \oplus &~&~\\
 ~ & ~ & |[name=a2]| 0& |[name=b2]|0 &~&~\\
        };

        \path[overlay,->, font=\scriptsize,>=latex]
        (ka) edge (kb)
         (kb) edge (kc)
        
          (kd) edge (ke)
          (ke) edge (kf)
          
            (kb) edge (a1)
            (kb) edge (a2)
             (a1) edge (b1)
 (a1) edge (kd)
 (a2) edge (b2)
(a2) edge (kd)
(kc) edge (b1)
(kc) edge (b2)
(b1) edge (ke)
(b2) edge (ke)
 ;
\end{tikzpicture}
\caption{$3$ negative exponents}\label{check2}
\end{figure}
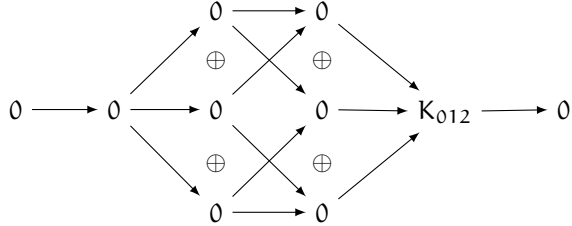
\item[$2$ negative exponents] The monomial $X_0^{a_0}\cdot X_1^{a_1}\cdot X_2^{a_2}$ where two exponents are negative, say $a_0,a_1<0$. Then we can perfectly lift to coboundary coming from $U_0\cap U_1$, which gives exactness.

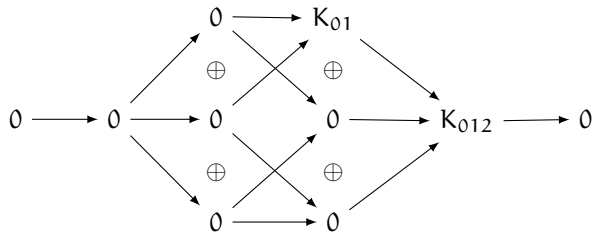
\begin{figure}[H]
\centering
 \begin{tikzpicture}
 []
        \matrix (m) [
            matrix of math nodes,
            row sep=0.5em,
            column sep=2.5em,
                   ]
{ ~ & ~ &  |[name=a1]|0 & |[name=b1]|K_{01} &~&~\\
 ~ & ~ & \oplus & \oplus &~&~\\
 |[name=ka]| 0 &|[name=kb]| 0 & |[name=kc]|0 &|[name=kd]|0 &|[name=ke]|K_{012}&|[name=kf]|0\\
  ~ & ~ & \oplus & \oplus &~&~\\
 ~ & ~ & |[name=a2]| 0& |[name=b2]|0 &~&~\\
        };

        \path[overlay,->, font=\scriptsize,>=latex]
        (ka) edge (kb)
         (kb) edge (kc)
        
          (kd) edge (ke)
          (ke) edge (kf)
          
            (kb) edge (a1)
            (kb) edge (a2)
             (a1) edge (b1)
 (a1) edge (kd)
 (a2) edge (b2)
(a2) edge (kd)
(kc) edge (b1)
(kc) edge (b2)
(b1) edge (ke)
(b2) edge (ke)
 ;
\end{tikzpicture}
\caption{$2$ negative exponents}\label{check3}
\end{figure}

\pagebreak

\item[$1$ negative exponent] The monomial $X_0^{a_0}\cdot X_1^{a_1}\cdot X_2^{a_2}$ where one exponents is negative, say $a_0<0$, we get the complex (Figure \ref{check4}). Notice that $K_0$ maps injectively giving zero cohomology group.

\begin{figure}[H]
\centering
 \begin{tikzpicture}
 []
        \matrix (m) [
            matrix of math nodes,
            row sep=0.5em,
            column sep=2.5em,
                   ]
{ ~ & ~ &  |[name=a1]|K_0 & |[name=b1]|K_{01} &~&~\\
 ~ & ~ & \oplus & \oplus &~&~\\
 |[name=ka]| 0 &|[name=kb]| 0 & |[name=kc]|0 &|[name=kd]|K_{02} &|[name=ke]|K_{012}&|[name=kf]|0\\
  ~ & ~ & \oplus & \oplus &~&~\\
 ~ & ~ & |[name=a2]| 0& |[name=b2]|0 &~&~\\
        };

        \path[overlay,->, font=\scriptsize,>=latex]
        (ka) edge (kb)
         (kb) edge (kc)
        
          (kd) edge (ke)
          (ke) edge (kf)
          
            (kb) edge (a1)
            (kb) edge (a2)
             (a1) edge (b1)
 (a1) edge (kd)
 (a2) edge (b2)
(a2) edge (kd)
(kc) edge (b1)
(kc) edge (b2)
(b1) edge (ke)
(b2) edge (ke)
 ;
\end{tikzpicture}
\caption{$1$ negative exponent}\label{check4}
\end{figure}
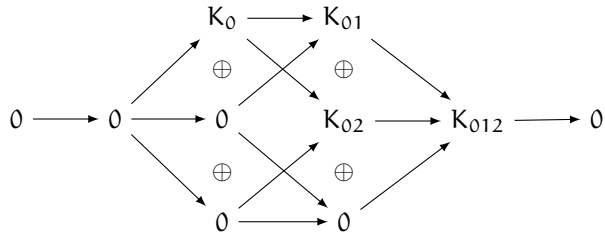

Furthermore, the mapping in the Figure \ref{check5} gives Kernel when $f=g$ which is possible for zero only. Again giving us zero cohomology groups.
\begin{figure}[H]
\centering
 \begin{tikzpicture}
 []
        \matrix (m) [
            matrix of math nodes,
            row sep=0.5em,
            column sep=2.5em,
                   ]
{ ~ & ~ &  |[name=a1]|K_0 & |[name=b1]|f &~&~\\
 ~ & ~ & \oplus & \oplus &~&~\\
 |[name=ka]| 0 &|[name=kb]| 0 & |[name=kc]|0 &|[name=kd]|g &|[name=ke]|f-g&|[name=kf]|0\\
  ~ & ~ & \oplus & \oplus &~&~\\
 ~ & ~ & |[name=a2]| 0& |[name=b2]|0 &~&~\\
        };

        \path[overlay,->, font=\scriptsize,>=latex]
        (ka) edge (kb)
         (kb) edge (kc)
        
          (kd) edge (ke)
          (ke) edge (kf)
          
            (kb) edge (a1)
            (kb) edge (a2)
             (a1) edge (b1)
 (a1) edge (kd)
 (a2) edge (b2)
(a2) edge (kd)
(kc) edge (b1)
(kc) edge (b2)
(b1) edge (ke)
(b2) edge (ke)
 ;
\end{tikzpicture}
\caption{Mapping for $1$ negative exponent}\label{check5}
\end{figure}
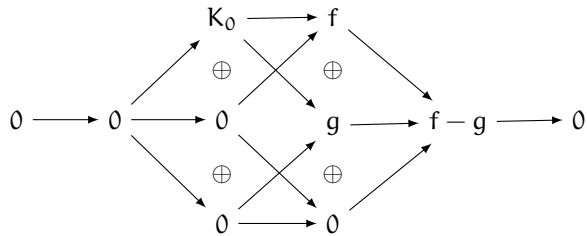

\pagebreak

\item[$0$ negative exponent] The monomial $X_0^{a_0}\cdot X_1^{a_1}\cdot X_2^{a_2}$ where none of the exponents is negative $a_i>0$, gives the complex Figure \ref{check7}.

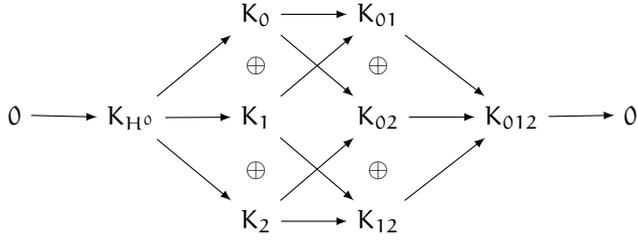
\begin{figure}[H]
\centering
 \begin{tikzpicture}
 []
        \matrix (m) [
            matrix of math nodes,
            row sep=0.5em,
            column sep=2.5em,
                   ]
{ ~ & ~ &  |[name=a1]|K_0 & |[name=b1]|K_{01} &~&~\\
 ~ & ~ & \oplus & \oplus &~&~\\
 |[name=ka]| 0 &|[name=kb]| K_{H^0} & |[name=kc]|K_1 &|[name=kd]|K_{02} &|[name=ke]|K_{012}&|[name=kf]|0\\
  ~ & ~ & \oplus & \oplus &~&~\\
 ~ & ~ & |[name=a2]| K_2& |[name=b2]|K_{12} &~&~\\
        };

        \path[overlay,->, font=\scriptsize,>=latex]
        (ka) edge (kb)
         (kb) edge (kc)
        
          (kd) edge (ke)
          (ke) edge (kf)
          
            (kb) edge (a1)
            (kb) edge (a2)
             (a1) edge (b1)
 (a1) edge (kd)
 (a2) edge (b2)
(a2) edge (kd)
(kc) edge (b1)
(kc) edge (b2)
(b1) edge (ke)
(b2) edge (ke)
 ;
\end{tikzpicture}
\caption{$0$ negative exponent}\label{check7}
\end{figure}

Consider the SES of complex as in Figure \ref{check8}  . The top and bottom row come from the $1$ negative exponent case, thus giving zero cohomology. The SES of complex gives LES of cohomology groups, since top and bottom row have zero cohomology, so does the middle.

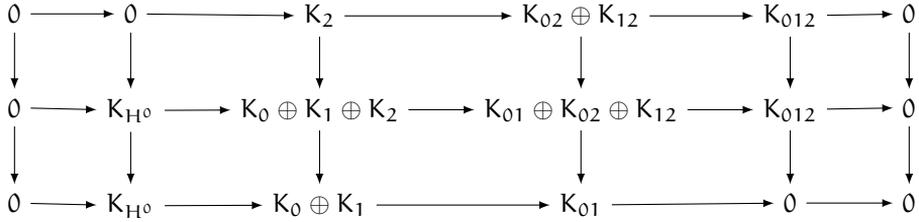
\begin{figure}[H]
\centering
 \begin{tikzpicture}\label{check8}
 []
        \matrix (m) [
            matrix of math nodes,
            row sep=2em,
            column sep=2.5em,
                   ]
{   |[name=aa]|0 &  |[name=ab]|0 &  |[name=a1]|K_2 & |[name=b1]|K_{02}\oplus K_{12} &  |[name=c1]|K_{012}& |[name=d1]|0\\
 |[name=ka]| 0 &|[name=kb]| K_{H^0} & |[name=kc]|K_0\oplus K_1\oplus K_2 &|[name=kd]|K_{01}\oplus K_{02} \oplus K_{12} &|[name=ke]|K_{012}&|[name=kf]|0\\
 |[name=qa]| 0 &  |[name=qb]|K_{H^0} & |[name=a2]| K_0\oplus K_1& |[name=b2]|K_{01} & |[name=c2]|0& |[name=d2]|0\\
        };

        \path[overlay,->, font=\scriptsize,>=latex]
(aa) edge (ab)
(ab) edge (a1)
(a1) edge (b1)
(b1) edge (c1)
(c1) edge (d1)
        (ka) edge (kb)
        (kb) edge (kc)
        (kc) edge (kd)
        (kd) edge (ke)
        (ke) edge (kf)
(qa) edge (qb)
(qb) edge (a2)
(a2) edge (b2)
(b2) edge (c2)
(c2) edge (d2)

(aa) edge (ka)
(ab) edge (kb)
(a1) edge (kc)
(b1) edge (kd)
(c1) edge (ke)
(d1) edge (kf)

        (ka) edge (qa)
        (kb) edge (qb)
        (kc) edge (a2)
        (kd) edge (b2)
        (ke) edge (c2)
	    (kf) edge (d2)

 ;
\end{tikzpicture}
\caption{SES of Complex}\label{check8}  
  \end{figure}
\end{description}

\section{Kunneth Formula} We can produce a complex for $\perfp\times \perfpa{m}$ by taking (completed) tensor product of the corresponding \v{C}ech complex associated with each space, and by the Theorem of Eilenberg-Zilber we get
{\small
\begin{equation}
H^i(\perfp\times \perfpa{m},\curly{O}(a,b))=\sum_{j=0}^iH^j(\perfp,\curly{O}(a))\otimes H^{i-j}(\perfpa{m},\curly{O}(b))~~a,b\in\ds{Z}[1/p]
\end{equation}
}
Furthermore, we can define a cup product following \cite[pp 194]{liu2002algebraic} to get a homomorphism
{\small
\begin{equation}
\smile:H^p(\perfp,\curly{O}(a)) \times H^q(\perfpa{m},\curly{O}(b))\ra H^{p+q}(\perfp\times \perfpa{m},\curly{O}(a,b))~~a,b\in\ds{Z}[1/p]
\end{equation}

}

\section{Braided Dimension} The dimensions of the cohomology groups denoted by $h^i$ are infinite and thus it is not possible to use them for any direct computations. In particular the Euler characteristic does not make sense. The situation can be rectified by introduction of a grading on these dimensions which will be called Braided dimension, but the degree of all the grades is same.

Let me start with an example, consider $H^0(\pro,\curly{O}_{\pro}(2))$ is generated by infinitely many monomials as explained in section \ref{coho1}. But, these monomials can be arranged in terms of the denominator being a power of $p$ (with risk of repetition, which can be avoided if the reader prefers it). The monomial $x^ay^b$ will be denoted by $(a,b)$ and the monomial $x^{a/p^m}y^{b/p^m}$ will also be denoted by $(a,b)$ but will lie in the row with power of $p$ as $m$. In the table below we fix $p=3$ and consider $\deg=2$. For example in the first row $(1,1)$ corresponds to $xy$, in the second row $(2,4)$ corresponds to $x^{2/3}y^{4/3}$, and the third row $(2,16)$ corresponds to $x^{2/9}y^{16/9}$. Notice the repetition in the second row where $(3,3)$ corresponds to $x^{3/3}y^{3/3}=xy$. This can be avoided by simply dropping the term, but has been kept to ease computations.

\begin{equation}\label{grd1}
\begin{aligned}
\text{power of }p & ~~\text{monomials}& ~~\text{dim}\\
0 \qquad& (2,0),(1,1),(0,2)& ~~3\\
1 \qquad& (0,6)(1,5),(2,4),(3,3),(4,2), (5,1), (6,0)& ~~7=(3\times 2)+1\\
2 \qquad& (0,18),(1,17),(2,16),(3,15),(4,14),\ldots & ~~19=(3^2\times 2)+1\\
\vdots \qquad& \qquad\qquad\vdots & ~~\vdots
\end{aligned} 
\end{equation}

The point I am trying to make is that for every grade the dimension is finite, in the first row partition of $3$ is considered, partition of $6$ in the second row, partition of $18$ in the third row and so on. We can write this as a tuple
\begin{equation}
(3,7,19,\ldots, (3^n\times\deg)+1,\ldots)\text{ here }\deg=2
\end{equation}

Still working in the space $\pro$ consider $\deg=s$ with fixed prime $p$, the tuple is
\begin{equation}\label{eq:pdim}
\left(\binom{1+s}{1},ps+1,p^2 s+1,\ldots, (p^n\times\deg)+1,\ldots\right)\text{ here }\deg=s
\end{equation}
and for the sake of computation one could add or subtract the above to say another tuple with $\deg =t$. The first co-ordinate of the tuple always recovers the case of integer degrees. This above recipe can be extended to $\perfp$ and any degree $d$ for a fixed prime $p$. For example, for $h^0(\perfp)$, the number of variables are $n+1$ corresponding to $x_0\ldots,x_n$ and we get the Braided dimension as, 

\begin{equation}
\begin{aligned}
\text{power of }p & ~~\text{dim}\\
0\qquad & ~~\binom{n+d}{n}\\
1\qquad & ~~\text{number of partition of $pd$ into tuple of size $n+1$}\\
2 \qquad& ~~\text{number of partition of $p^2d$ into tuple of size $n+1$}\\
\vdots\qquad & \qquad\qquad\vdots 
\end{aligned} 
\end{equation}

All these partitions are finite in number and similar situation holds for $h^n$. Hence, the Euler characteristic now makes sense as a tuple. It is just addition/subtraction of tuples and the first co-ordinate of the tuple always recovers the case of standard Euler Characteristic.

\subsection{Braided dimension for fractions} Since, $\perfp$ is richer than its algebraic geometry counterpart, it affords rational degrees of the form $m/p^i$. It is possible to put grading on these, by simply considering the grading of integer $m$ and then multiplying that grade by $1/p^i$. For example consider $\deg =2/3$ and compare with \eqref{grd1}.
{\small
\begin{equation}
\begin{aligned}
\text{power of }p & ~~\text{monomials}& ~~\text{dim}\\
1\qquad &(2/3,0),(1/3,1/3),(0,2/3)& ~~3\\
2\qquad & (0,6/9)(1/9,5/9),(2/9,4/9),(3/9,3/9),(4/9,2/9), (5/9,1/9), (6/9,0)& ~~7=(3\times 2)+1\\
3\qquad & (0,18/27),(1/27,17/27),(2/27,16/27),(3/27,15/27),(4/27,14/27),\ldots & ~~19=(3^2\times 2)+1\\
\vdots\qquad & ~~~\qquad\vdots & ~~\vdots
\end{aligned} 
\end{equation}
}

Notice that the dimensions have not changed from the integer case. Following the convention of Braided degree defined before the fraction $1/3^i$ should have been left from the tuples, but is kept for clarity.
\subsection{Bezout Type computations}\label{bezoutexam} Look at the following computation as an example for a fixed prime $p=3$
\begin{equation}
h^1\curly{O}_{\pro}(-5)-h^1\curly{O}_{\pro}(-2)-h^1\curly{O}_{\pro}(-3)
\end{equation}

Table for $h^1\curly{O}_{\pro}(-5)$
\begin{equation}
\begin{aligned}
\text{power of }p & ~~\text{monomials}& ~~\text{dim}\\
0 \qquad& (-4,-1),(-3,-2), (-2,-3),(-1,-4),& ~~4\\
1\qquad & (-1,-14),(-2,-13),(-3,-12),(-4,-11), \ldots & ~~14=(3\times 5)-1\\
2 \qquad& (-1,-44),(-2,-43),(-3,-42),(-4,-41),\ldots & ~~44=(3^2\times 5)-1\\
\vdots\qquad &\qquad\qquad \vdots & ~~\vdots
\end{aligned} 
\end{equation}

Table for $h^1\curly{O}_{\pro}(-2)$
\begin{equation}
\begin{aligned}
\text{power of }p & ~~\text{monomials}& ~~\text{dim}\\
0\qquad & (-1,-1)& ~~1\\
1\qquad & (-1,-5),(-2,-4),(-3,-3),(-4,-2), (-5,-1) & ~~5=(3\times 2)-1\\
2\qquad & (-1,-17),(-2,-16),(-3,-15),(-4,-14),\ldots & ~~17=(3^2\times 2)-1\\
\vdots\qquad & \qquad\qquad\vdots & ~~\vdots
\end{aligned} 
\end{equation}

Table for $h^1\curly{O}_{\pro}(-3)$
\begin{equation}
\begin{aligned}
\text{power of }p & ~~\text{monomials}& ~~\text{dim}\\
0\qquad & (-1,-2), (-2,-1)& ~~2\\
1\qquad & (-1,-8),(-2,-7),(-3,-6),(-4,-5), \ldots & ~~8=(3\times 3)-1\\
2\qquad & (-1,-26),(-2,-25),(-3,-23),(-4,-22),\ldots & ~~26=(3^2\times 3)-1\\
\vdots & ~~~\vdots & ~~\vdots
\end{aligned} 
\end{equation}
The computation $h^1\curly{O}_{\pro}(-5)-h^1\curly{O}_{\pro}(-2)-h^1\curly{O}_{\pro}(-3)$ gives 
\begin{equation}
\begin{aligned}
&(4,14,44,\ldots)\\
-&(1,5,17,\ldots)\\
-&(2,8,26,\ldots)\\
=&(1,1,1,\ldots)\\
\end{aligned}
\end{equation}
This nothing short of a miracle and this would in fact hold for any equation of the form 
\begin{equation}
h^1\curly{O}_{\pro}(-s-t)-h^1\curly{O}_{\pro}(-s)-h^1\curly{O}_{\pro}(-t)
\end{equation}
because of the dimension computation given in \eqref{eq:pdim}. Further, computations of this sort can be carried for other degrees and for any $n$ in $\perfp$ and the first co-ordinate of the tuple is always known from algebraic geometry.
\subsubsection{Fraction Degree}

But more is possible, one could consider fractional degrees of the form $m/p^i$ and above would still hold, since the dimensions depend upon the numerator and not the denominator. For example, consider the degrees $-5/3,-2/3,-3/3$ for $p=3$ and the above computations hold, again a marvel which is not conceivable in algebraic geometry. 

Notice the fraction $-3/3$ or simply $-1$, this would have no sections in $\ds{P}^1$ but it does have sections in $\pro$. More concretely it is given as
{\small
\begin{equation}
\begin{aligned}
\text{power of }p & ~~\text{monomials}& ~~\text{sections}\\
1\qquad & (-1,-2), (-2,-1)& \frac{1}{x^{1/3}y^{2/3}},\frac{1}{x^{2/3}y^{1/3}},\ldots\\
2\qquad & (-1,-8),(-2,-7),(-3,-6), \ldots & ~~\frac{1}{x^{1/9}y^{8/9}},\frac{1}{x^{2/9}y^{7/9}},\frac{1}{x^{3/9}y^{6/9}},\ldots\\
3\qquad & (-1,-26),(-2,-25),(-3,-23),\ldots & ~~\frac{1}{x^{1/27}y^{26/27}},\frac{1}{x^{2/27}y^{25/27}},\frac{1}{x^{3/27}y^{24/27}},\ldots\\
\vdots\qquad & \qquad\qquad\vdots & ~~\vdots
\end{aligned} 
\end{equation}
}

\subsubsection{Vector bundles on $\pro$} Grothendieck's theorem on classification of vector bundles on $\ds{P}^1$ is a consequence of the fact that an integer can be partitioned only in finitely many ways. This no longer remains true in the fractional case, for example $2=\sum_{i\geq 0}2^{-i}$ for $p=2$, thus an infinite splitting is possible for every integer degree. It is not possible for every vector bundle to split over $\pro$ as shown in \cite{Kiran_AWS}, but using the concept of braided dimension, we can talk about splitting in every grade. The first grade then recovers Grothnedieck's classification.

\subsection{Euler Characteristic in Braided dimension} Recall that for a short exact sequence of sheaves 

\begin{equation}
0\ra \mathcal{F}\ra \mathcal{G}\ra\mathcal{H}\ra 0
\end{equation}
the identity $\chi(\mathcal{G})=\chi(\mathcal{F})+\chi(\mathcal{H})$ holds. In light of the example in section  \ref{bezoutexam} it is obvious to ask if this will hold for the entire Braided dimension.  Fortunately, it does hold justifying the idea of Braided dimension.

\subsubsection{Grading on Sheaves}
In Galois Theory we construct the following sequence of extensions
\begin{equation}
K[X]\subset K[X^{1/p}]\subset K[X^{1/p^2}]\subset \ldots\subset K[X^{1/p^i}]\subset\ldots
\end{equation}
with Noetherian ring every where except in the limit. This can be thought as a tuple
\begin{equation}\label{gal2}
\lr{K[X],K[X^{1/p}],K[X^{1/p^2}],\ldots,K[X^{1/p^i}],\ldots}.
\end{equation}
The above idea can now be carried to sheaves with notation $\ten{i}{\mathcal{F}}{}{}$ corresponding to the sheaf which would contain the $1/p^i$th power. For example, the tuple \eqref{gal2} can be thought of living on the perfectoid affine line as a Braided dimension sheaf as below.

\begin{equation}
\lr{K\lr{X},K\lr{X^{1/p}},K\lr{X^{1/p^2}},\ldots,K\lr{X^{1/p^i}},\ldots}.
\end{equation}
Hence, in Braided dimension the sheaf looks like the following:
\begin{equation}
\lr{\mathcal{F},~^1\mathcal{F},~^2\mathcal{F},\ldots,~^i\mathcal{F},\ldots}
\end{equation}
and the short exact sequence of sheaves in perfectoid spaces
\begin{equation}
0\ra\mathcal{F}\ra \mathcal{G}\ra \mathcal{H}\ra 0
\end{equation}
should be thought of as the following (with the first tuple corresponding to integer degrees)
\begin{figure}[H]
\centering
 \begin{tikzpicture}
 []
        \matrix (m) [
            matrix of math nodes,
            row sep=1.5em,
            column sep=2.5em,
                   ]
{
\langle &|[name=k1]|0 & |[name=k2]|0 &  |[name=k3]|0 & |[name=k4]|\cdots &|[name=k5]|0&|[name=k6]|\cdots & \rangle\\ 
\langle &|[name=a1]|\mathcal{F} & |[name=a2]|~^1\mathcal{F} &  |[name=a3]|~^2\mathcal{F} & |[name=a4]|\cdots &|[name=a5]|~^i\mathcal{F}&|[name=a6]|\cdots & \rangle\\
\langle &|[name=b1]|\mathcal{G} & |[name=b2]|~^1\mathcal{G} &  |[name=b3]|~^2\mathcal{G} & |[name=b4]|\cdots &|[name=b5]|~^i\mathcal{G}&|[name=b6]|\cdots& \rangle\\
\langle &|[name=c1]|\mathcal{H} & |[name=c2]|~^1\mathcal{H} &  |[name=c3]|~^2\mathcal{H} & |[name=c4]|\cdots &|[name=c5]|~^i\mathcal{H}&|[name=c6]|\cdots& \rangle\\
\langle &|[name=d1]|0 & |[name=d2]|0 &  |[name=d3]|0 & |[name=d4]|\cdots &|[name=d5]|0&|[name=d6]|\cdots& \rangle\\
        };

\path[overlay,->, font=\scriptsize,>=latex]
(k1) edge (a1)
(k2) edge (a2)
(k3) edge (a3)
(k5) edge (a5)      

(a1) edge (b1)
(a2) edge (b2)
(a3) edge (b3)
(a5) edge (b5)

(b1) edge (c1)
(b2) edge (c2)
(b3) edge (c3)
(b5) edge (c5)

(c1) edge (d1)
(c2) edge (d2)
(c3) edge (d3)
(c5) edge (d5)          
 ;
\end{tikzpicture}
\caption{Short Exact Sequence in Braided Dimension}\label{gradedses1}
\end{figure}

The figure \ref{gradedses1} makes it obvious that identity $\chi(~^i\mathcal{G})=\chi(~^i\mathcal{F})+\chi(~^i\mathcal{H})$ will hold in every grade. The perfectoid Euler characteristic can now be defined as the limit of individual euler characters if it exists. Notice that the Euler Characteristic might be infinite, as already shown for the $\perfp$ case, but the identity $\chi(~^i\mathcal{G})=\chi(~^i\mathcal{F})+\chi(~^i\mathcal{H})$ still might hold as $i\ra \infty$ giving meaningful computations.

\section{Bezout Theorem in $\perfpa{2}$} This section follows Miles Reid lectures in \cite[pp 32]{kollar2000complex}. For any two given polynomials of degree $m$ and $n$ respectively and high enough $d$ we have the sequence of vector spaces 

\begin{equation}
0\ra V_{d-n-m}\xra{-g,f}V_{d-n}\oplus V_{d-m}\xra{\begin{matrix}
f\\g
\end{matrix}}V_d\ra\oplus_P\curly{O}_P/\lr{f,g}\ra 0
\end{equation}
where $V_d=K\lr{X^{1/p^\infty},Y^{1/p^\infty},Z^{1/p^\infty}}_d$ (the $d$ th graded part  as in Reid). These vector spaces are finite dimensional in algebraic geometry and the euler characteristic computation gives rise to Bezout formula. But, these spaces are infinite dimensional in perfectoid case and thus the concept of braided dimension comes to rescue; there is a finite dimension for every grade and the first grade recovers $nm$, the bezout formulae.

\subsection{Partition of Integers} The general philosophy or slogan of braided dimension is `integer partitions determine geometry'. Given a degree $d$ we can partition it in $d+1$ two co-ordinates as 

\begin{equation}
(d,0),(d-1,1),(d-2,2), \ldots, (1,d-1),(0,d)
\end{equation}

The partition into three tuples would require the second co-oordinate to be a constant, and partition the first co-ordinate into two co-ordinates, that is partition $(d,0)$ as 

\begin{equation}
(d,0,0),(d-1,1,0),(d-2,2,0), \ldots, (1,d-1,0),(0,d,0)
\end{equation}
a, total of $d+1$ and so on.

Thus, for three co-ordinates below (each for power of $X,Y$ and $Z$) the total matches with $\binom{d+2}{2}$ the dimension of $V_d$ in algebraic geometry.

\begin{equation}
\begin{aligned}
\text{Two Co-ords}\qquad &\text{number of $3$ co-ords}\\
(d,0)\qquad\qquad & \qquad d+1\\
(d-1,1)\qquad\qquad & \qquad d\\
(d-2,2)\qquad\qquad & \qquad d-1\\
\vdots\qquad\qquad & \qquad \vdots\\
(0,d)\qquad\qquad & \qquad 1\\
\text{Total}\qquad &  \frac{(d+1)(d+2)}{2}=\binom{d+2}{2}\\
\end{aligned}
\end{equation}

In the perfectoid case in the grade corresponding to $1/p$ the degree considered is $dp$, in the next grade it is $dp^2$ and continue. The above recipe can be used to compute Euler Character in every grade and thus the dimension. 

\begin{equation}
\begin{aligned}
\dim_K\bigoplus_P\curly{O}_P/\lr{f,g}&=\binom{d-n-m}{2}-\binom{d-n}{2}-\binom{d-m}{2}+\binom{d}{2}\\
\dim_K~^1\bigoplus_P\curly{O}_P/\lr{f,g}&=\binom{p(d-n-m)}{2}-\binom{p(d-n)}{2}-\binom{p(d-m)}{2}+\binom{pd}{2}\\
\dim_K~^2\bigoplus_P\curly{O}_P/\lr{f,g}&=\binom{p^2(d-n-m)}{2}-\binom{p^2(d-n)}{2}-\binom{p^2(d-m)}{2}+\binom{p^2d}{2}\\
\vdots\qquad&=\qquad\vdots
\end{aligned}
\end{equation}

The first line is algebraic geometry, and the rest is perfectoid version, which depends upon integer partitions as given in \cite{ram1}.

\subsection{Perfectoid Veronese Embedding}
Consider the $2$ uple embedding given by 
\begin{equation}
[x:y]\mapsto [x^2:xy:\ldots:x^{a/p^i}y^{a/p^i}:\ldots y^2] \text{ where }\frac{a}{p^i}+\frac{b}{p^i}=2
\end{equation}
The map is from $\pro\ra\perfpa{\infty}$, this does not seem to much useful. The previous section gives a hint to define this map in various pieces depending upon the power of the prime $p$.

Fix $p=3$ and consider the following maps based upon partitions given in previous section

\begin{equation}
\begin{aligned}
[x:y]\mapsto & [x^2:xy:y^2]\\
&[x^2:x^{1/3}y^{1/3}:x^{2/3}y^{4/3}:xy:x^{4/3}y^{2/3}: x^{5/3}y^{1/3}:y^2]\\
&[x^2:x^{1/9}y^{17/9}:x^{2/9}y^{16/9}:\ldots:y^2]\\
\end{aligned}
\end{equation}

Thus there is a map $\pro\ra\perfpa{2},\perfpa{6}, \perfpa{18},\ldots,\perfpa{3^n\times 2},\ldots$. Notice that the inclusion $\perfpa{2}\subset \perfpa{6}\subset \perfpa{18}\ldots$. This toy example now clearly generalizes to the standard $d$ uple embedding in $n+1$ variables and $N$ monomials where $N=\binom{n+d}{n}-1$. The mapping will $\perfpa{n}\ra\perfpa{N},\perfpa{N_p},\perfpa{N_{p^2}},\ldots  $ where $\perfpa{N_{p^i}} $ means the embedding obtained by considering partitions of degree as done in previous section and the example above. The above construction can be used to put a space inside a space inside a space thus constructing an infinite tower of spaces, where each space is dependent upon the one before it.

\subsubsection{Morphisms into $\perfp$} 

In the Braided dimension system the invertible sheaf $\curly{O}_{\perfp}(1)$ is generated in different grades (although it has infinitely many sections). For example for $\pro$ and $p=3$ the Braided sections of $\curly{O}_{\pro}(1)$ are given as
\begin{equation}\label{amp2}
[x_0:x_1][x_0:x_0^{1/3}x_1^{2/3}:x_0^{2/3}x_1^{1/3}:x_1][x_0:x_0^{1/9}x_1^{8/9}:\ldots:x_1]\ldots
\end{equation}

Now let $\phi:X\ra\perfp$ be a $K$ morphism of $X$ into $\perfp$, then we can construct a Braided invertible sheaf $\curly{L}$ on $X$ given as $\phi^*(\curly{O}(1))$ with sections $s_i=\phi^*(x_i)$, again the $s_i$ should be considered living in various grades.

It is also possible to consider rational degrees on $\perfp$, for example, $\curly{O}_{\pro}(1/3)$ on $\pro$ would just be
\begin{equation}
\begin{aligned}
[x_0^{1/3}:x_1^{1/3}][x_0^{1/3}:x_0^{1/9}x_1^{2/9}:x_0^{2/9}x_1^{1/9}:x_1^{1/3}][x_0^{1/3}:x_0^{1/27}x_1^{8/27}:\ldots:x_1^{1/3}]\ldots
\end{aligned}
\end{equation}

In the new brave world of $\perfp$ one could take Braided invertible sheaf $\phi^*(\curly{O}(1/3))$, which is not possible in standard algebraic geometry.

\subsubsection{Ample Vector Bundles}
It is evident that $\curly{O}_{\perfp}(1)$ is ample since it is generated by the global sections and one could pull back these line bundles to a perfectoid space $X\ra\perfp$, that is $\phi^*\curly{O}(1)$ as given on \cite[pp 150]{hartshorne1977algebraic}). The concept of braided dimension gives a precise structure to this pull back. For example for $\perfpa{1}$ we are pulling back the tower of infinite spaces as given in \eqref{amp2}. The concept of braided dimension helps transfer the constructions directly from algebraic geometry into perfectoid geometry for every grade.

But very ample requires embedding into a projective space, which cannot be done in $\perfp$ for $n$ finite, that is we no longer have a embedding defined by hyperplanes $[x_0:\ldots:x_n]$. One way to fix it would be to simply take a finite number of sections and define an embedding, but that would be arbitrary. Thus, a more useful notion for very ample is required.

\begin{mdef}
$\curly{O}_{\perfp}(d)$ is very ample if it can be embedded into infinitely many copies of projective space as a Perfectoid Veronese Embedding. 
\end{mdef}

Notice that $\curly{O}_{\perfp}(d)$ is ample iff it is very ample iff $d>0$.

The  very important notion of ampleness has been rescued and one can now look at all the ampleness theorem analogues in algebraic geometry, apply the tilting functor and see what the result would be over a finite field.

\section{Intersection Multiplicities in perfectoid affine plane } Another application of Braided dimension is in intersection theory. Of course, there is a well defined notion of a local ring at point $P$ in perfectoid affine plane $\arov{2}$. 

\begin{equation}
\curly{O}_P:=\left\{\frac{f}{g}\text{ such that }g(P)\neq 0 \text{ and  }f,g\in K\lr{X^{1/p^\infty},Y^{1/p^\infty}}\right\}
\end{equation}

The evaluation map at point $P$ gives us a ring homomorphism $\curly{O}_P\rightarrow K$. The kernel of this map is the maximal ideal of the ring which will be denoted by $\id{m}_P$. 

The ring $\curly{O}_P$ can be rewritten in the Braided dimension form with the first tuple corresponding to $p=0$ and the the $i$th tuple corresponding to $p=i$. More formally 
\begin{equation}
\left(\ten{}{\curly{O}}{P}{}, \ten{1}{\curly{O}}{P}{},\ten{2}{\curly{O}}{P}{}\ldots,\ten{i}{\curly{O}}{P}{},\ldots\right)
\end{equation}

For a given curve $F=0$ (with integer powers), we have to consider all the other curves that come with it $F^{1/p}=0,\ldots, F^{1/p^i}=0,\ldots$. The reader could keep the perfectoid field $\ds{F}_p((t^{1/p^\infty}))$ in mind as a crutch for computations in $p$th power roots. The newton binomial formula for $p$th roots should be used carefully, since it can give series whose coefficients might not converge to zero. Alternatively one could instead define $F^{1/p^i}:=F(X^{1/p^i},Y^{1/p^i})$, which would give a simpler but a different notion and consider the corresponding ideal generated as $\lr{F(X^{1/p^i},Y^{1/p^{i-1}}),F(X^{1/p^{i-1}},Y^{1/p^{i-1}}),\ldots,F(X,Y)}$.

\begin{mdef}
Intersection multiplicity at a point $P$ for curves $F$ and $G$ (with integer powers) is defined as a tuple
\begin{equation}
\mu_P(F,G):=\left\{\dim_K \frac{\curly{O}_{P}}{\lr{F,G}},\dim_K \frac{\ten{1}{\curly{O}}{P}{}}{\lr{I^{1/p}}}, \ldots,\dim_K \frac{\ten{i}{\curly{O}}{P}{}}{\lr{I^{1/p^i}}},\ldots  \right\}
\end{equation}

where ${\ten{1}{\curly{O}}{P}{}}/{\lr{I^{1/p}}}$ stands for four rings (with risk of confusion)\\
\begin{equation}
\frac{\ten{1}{\curly{O}}{P}{}}{\lr{F^{1/p},G^{1/p}}},\frac{\ten{1}{\curly{O}}{P}{}}{\lr{F^{1/p},G}},\frac{\ten{1}{\curly{O}}{P}{}}{\lr{F,G^{1/p}}},\frac{\ten{1}{\curly{O}}{P}{}}{\lr{F,G}}
\end{equation}

or more generally ${\ten{i}{\curly{O}}{P}{}}/{\lr{I^{1/p^i}}}$ stands for the tuple of following rings\\
\begin{equation}
\frac{\ten{i}{\curly{O}}{P}{}}{\lr{F^{1/p^i},G^{1/p^i}}},\frac{\ten{1}{\curly{O}}{P}{}}{\lr{F^{1/p^i},G^{1/p^{i-1}}}},\ldots,\frac{\ten{1}{\curly{O}}{P}{}}{\lr{F,G}}
\end{equation}

The tuples are arranged in an order starting with lowest power of $F$ and then considering all the corresponding powers of $G$.

\end{mdef}

In case $F$ and $G$ already have fractional powers, begin with the tuple where they belong and if possible compute their $p^i$th power to fill in the tuples before it, otherwise set previous entries as zero.

An immediate consequence of the above is that $\mu_P(F, G)=\mu_P(G, F)$ only in the integer degree grade. The ordering sequence between $F$ and $G$ will have to be switched for the equality to hold in other grades.

Similarly, the changes in ordering suggested above lead the equality $\mu_P( F, G)=\mu_P( F, G+FH)$ (follows from $\lr{F,G}=\lr{F,G+FH}$). Furthermore if $F,G$ have a common component then $\mu_P(F,G)=\infty$.

\begin{exam}[Intersection Multiplicity of coordinate axis]
Consider the evaluation map $\curly{O}_0\ra K$ the kernel of this map is precisely the fractions of the form $f/g$ such that $f(0)=0$ (that is no constant term). 

\begin{equation}
\mu_0(X,Y):=\left\{\dim_K \frac{K\lr{X,Y}}{\lr{X,Y}},\dim_K \frac{K\lr{X^{1/p}Y^{1/p}}}{\lr{I^{1/p}}}, \ldots,\dim_K \frac{K\lr{X^{1/p^i}Y^{1/p^i}}}{\lr{I^{1/p^i}}},\ldots  \right\}
\end{equation}

where $\dim_K{K\lr{X^{1/p}Y^{1/p}}}/{\lr{I^{1/p}}}$ stands for four rings \\
\begin{equation}
\dim_K\frac{K\lr{X^{1/p}Y^{1/p}}}{\lr{X^{1/p},Y^{1/p}}},\dim_K\frac{K\lr{X^{1/p}Y^{1/p}}}{\lr{X^{1/p},Y}},\dim_K\frac{K\lr{X^{1/p}Y^{1/p}}}{\lr{X,Y^{1/p}}},\dim_K\frac{K\lr{X^{1/p}Y^{1/p}}}{\lr{X,Y}}=(1,p,p,p^2)
\end{equation}

This gives us the dimension as $(1,(1,p,p,p^2),\ldots)$, notice the four enteries in the second co-ordinate. A physicist could now describe a web of events all of which start at $1$. Of course the special interest is in the sequence in which we extract $1$ to get $(1,1,1,\ldots,)$ which comes from 

\begin{equation}\label{mult1}
\left\{\dim_K \frac{K\lr{X,Y}}{\lr{X,Y}},\dim_K \frac{K\lr{X^{1/p}Y^{1/p}}}{\lr{X^{1/p},Y^{1/p}}}, \ldots,\dim_K \frac{K\lr{X^{1/p^i}Y^{1/p^i}}}{\lr{X^{1/p^i},Y^{1/p^i}}},\ldots  \right\}
\end{equation} 

\end{exam}

In light on \eqref{mult1} special perfectoid multiplicity is defined below.

\begin{mdef}
The special perfectoid intersection multiplicity at a point $P$ is defined as a limit (if it exists) of the following sequence
\begin{equation}
\mu_P(F,G):=\left\{\dim_K \frac{\curly{O}_{P}}{\lr{F,G}},\dim_K \frac{\ten{1}{\curly{O}}{P}{}}{\lr{F^{1/p},G^{1/p}}}, \ldots,\dim_K \frac{\ten{i}{\curly{O}}{P}{}}{\lr{F^{1/p^i},G^{1/p^i}}},\ldots  \right\}
\end{equation}

\end{mdef}

The above definition puts a handle on the concept of intersection multiplicity and helps interpret theorems from algebraic geometry in perfectoid geometry.

\section{Perfect Pairing} Recall the perfect pairing given on \cite[pp 225]{hartshorne1977algebraic}, Theorem 5.1(d) for $X=\ds{P}_A^r$.

\begin{equation}
H^0(X,\curly{O}_X(n))\times H^r(X,\curly{O}_X(-n-r-1)) \ra H^r(X,\curly{O}_X(-r-1))\cong A
\end{equation}

It was already shown that this pairing does not hold for $\perfp$. In this section we investigate this under the lens of Braided Dimension and realize that perfect pairing is a mirage which disappears when looked at deeper layers.

For example setting $n=2$ and $r=2$ we have the perfect pairing

\begin{align}\label{pair1}
\{x^2,y^2,z^2,xy,yz,xz\}\times \left\{\frac{1}{x^3yz},\frac{1}{xy^3z},\frac{1}{xyz^3},\frac{1}{x^2y^2z},\frac{1}{xy^2z^2},\frac{1}{x^2yz^2}\right\}=\frac{1}{xyz}
\end{align}
and such a pairing will always hold for an integer grade. In particular if $s$ is in the first set on LHS above, then its counterpart in the second set is $1/sxyz$. This pairing relies on the fact that there is only one element $1/xyz$ on the right hand side. 

In the perfectoid case say for $p=3$, the RHS of \eqref{pair1} consists of non zero partitions of $\deg\cdot p$ which in our example is $3\cdot 3=9$, and some of the elements in grade corresponding to $1/p$ are listed below
\begin{equation}
\begin{aligned}
\frac{1}{xyz} &\longleftrightarrow (-3,-3,-3)\\
\frac{1}{x^{1/3}y^{1/3}z^{1/3}z^2} &\longleftrightarrow (-1,-1,-7)\\
\frac{1}{x^{1/3}y^{2/3}z^2} &\longleftrightarrow (-1,-2,-6)\\
\vdots\qquad & \qquad\vdots
\end{aligned}
\end{equation}
Denoting the element of the base as $b$ (RHS of the column) and global section as $s$ there is a paring for every $\{s, 1/sb\}$ which gives $b$, although we have multiple elements of $b$. For example setting $s=x^2$ the following is obtained 

\begin{equation}\label{pair2}
\begin{aligned}
s\qquad&\frac{1}{sb}&\longleftrightarrow& \qquad b\\
x^2\qquad&\frac{1}{x^3yz} &\longleftrightarrow& \qquad\frac{1}{xyz} \\
x^2\qquad&\frac{1}{x^2x^{1/3}y^{1/3}z^{1/3}z^2} &\qquad\longleftrightarrow& \frac{1}{x^{1/3}y^{1/3}z^{1/3}z^2}\\
x^2\qquad&\frac{1}{x^2x^{1/3}y^{2/3}z^2} &\longleftrightarrow& \qquad\frac{1}{x^{1/3}y^{2/3}z^2}\\
\vdots\qquad & \qquad\vdots
\end{aligned}
\end{equation}

In the integer degree case $b$ was fixed as $1/xyz$ and the pairing $\{s, 1/sb\}$ was obtained by varying the section $s$, giving a `perfect pairing'. In the perfectoid case as shown above in \eqref{pair2}, we can fix $s$ and vary the base $b$ and still get a perfectly valid pairing. Thus, in the example stated above perfect pairing is consequence of the fact that $\deg 2-\deg 5=-\deg 3$ (nothing deeper).

\section{Blow Up} We start with a basic examples of blowing up a point on perfectoid varieties. These examples make sense only in a perfectoid world and not as a polynomial (variety in Algebraic Geometry).

\subsubsection{Example}
Let the variety $V$ be $y^{1/4}=x^{1/4}-x^{1/2}$ over perfectoid field (say $\rwhat{\ds{Q}_2(\mu_2^\infty)}$). We will blow up this variety over $(0,0)$.

\begin{equation}
\{(x,y)\times(u:v)\text{ such that } xv=uy,y^{1/4}=x^{1/4}-x^{1/2},(x,y)\neq(0,0)\}
\end{equation}

We look at the cases $u=1$ and $v=1$. We start by considering the first case $u=1$, which gives us $xv=y$ which can be plugged into $y^{1/4}=x^{1/4}-x^{1/2}$ to get 

\begin{equation}
\begin{aligned}
x^{1/4}v^{1/4}&=x^{1/4}-x^{1/2}\\
v^{1/4}&=1-x^{1/4}\\
v^{1/4}\ra 1 &\text{ as }x\ra 0
\end{aligned}
\end{equation}

The blow up over the origin are the points $(0,0)\times (1:\alpha)$ where $\alpha$ is solution to $v^{1/4}=1$.

We now consider the case $v=1$, which gives $x=yu$ which can be plugged into $y^{1/4}=x^{1/4}-x^{1/2}$ to get 
\begin{equation}
\begin{aligned}
y^{1/4}&=y^{1/4}u^{1/4}-y^{1/2}u^{1/2}\\
1&=u^{1/4}-y^{1/4}u^{1/2}\\
u^{1/4}\ra 1 &\text{ as }y\ra 0
\end{aligned}
\end{equation}
The blow up over the origin are the points $(0,0)\times (\beta:1)$ where $\beta$ is solution to $u^{1/4}=1$.

\subsubsection{Example} We now come to the standard example $y^2=x^3$ which we blow up at $(0,0)$. We consider two forms of this $y=x^{3/2}$ over $\rwhat{\ds{Q}_2(\mu_2^\infty)}$ and  $y^{2/3}=x$ in $\rwhat{\ds{Q}_3(\mu_3^\infty)}$

\begin{equation}
\{(x,y)\times(u:v)\text{ such that } xv=uy,y=x^{3/2},(x,y)\neq(0,0)\}
\end{equation} 

For $u=1$ we get $xv=y$ which plugged into $y=x^{3/2}$ gives $xv=x^{3/2}$ or $v=x^{1/2}$, thus as $x\ra 0$ we get $v\ra 0$. Hence the blow up over origin is $(0,0)\times (1:0)$. (Same as in case $y^2=x^3$).

For $v=1$ we get $x=yu$ which plugged into $y=x^{3/2}$ gives $y=y^{3/2}u^{3/2}$ or $1=y^{1/2}u^{3/2}$. No value of $u$ would satisfy the above if $y\ra 0$. Hence, for $v=1$ there are no points over the origin. (Same as in case $y^2=x^3$).

We now consider the other curve 

\begin{equation}
\{(x,y)\times(u:v)\text{ such that } xv=uy,y^{2/3}=x,(x,y)\neq(0,0)\}
\end{equation} 

For $u=1$ we get $xv=y$ which plugged into $y^{2/3}=x$ gives $x^{2/3}v^{2/3}=x$ or $v^{2/3}=x^{1/3}$, thus as $x\ra 0$ we get $v\ra 0$. Hence the blow up over origin is $(0,0)\times (1:0)$. (Same as in case $y^2=x^3$).

For $v=1$ we get $x=yu$ which plugged into $y^{2/3}=x$ gives $y^{2/3}=yu$ or $1=y^{1/3}u$. No value of $u$ would satisfy the above if $y\ra 0$. Hence, for $v=1$ there are no points over the origin. (Same as in case $y^2=x^3$).

\subsection{Blow Up Affine Plane} We describe blow up of a perfectoid affine plane over origin following \cite[p 163]{eisenbud2006geometry}. $\arov{2}=K\lr{X^{1/p^\infty}, Y^{1/p^\infty}}$, it is the union of two open sets $U'=K\lr{X_1^{1/p^\infty},Y_1^{1/p^\infty}}$ and $U''=K\lr{X_2^{1/p^\infty},Y_2^{1/p^\infty}}$ and consider the maps $\phi_1:U'\ra \arov{2}$ and $\phi_2:U''\ra\arov{2}$ which give us the ring homomorphisms.
\begin{equation}
\begin{aligned}
\phi_1^{\#}:K\lr{X^{1/p^\infty}, Y^{1/p^\infty}}&\ra K\lr{X_1^{1/p^\infty},Y_1^{1/p^\infty}}\\
X\mapsto X_1 &\text{ and } Y\mapsto X_1Y_1\\
K\lr{X^{1/p^\infty}, Y^{1/p^\infty},X^{-1/p^\infty},}&\ra K\lr{X_1^{1/p^\infty},Y_1^{1/p^\infty},X_1^{-1/p^\infty}}\label{eq:openset1}\\
X\mapsto X_1, & Y\mapsto X_1Y_1 \text{ and }YX^{-1}\mapsto Y_1
\end{aligned}
\end{equation}
From \eqref{eq:openset1}  we have an isomorphism between open subsets $U_{x}=K\lr{X^{1/p^\infty}, Y^{1/p^\infty},X^{-1/p^\infty}}$ and $U'_{x}=K\lr{X_1^{1/p^\infty},Y_1^{1/p^\infty},X_1^{-1/p^\infty}}$. Similarly, we have an isomorphism between open sets $U_y=K\lr{X^{1/p^\infty}, Y^{1/p^\infty},Y^{-1/p^\infty}} $ and $U'_{y}=K\lr{X_1^{1/p^\infty},Y_1^{1/p^\infty},Y_1^{-1/p^\infty}}$. In particular we have
\begin{equation}
U'_{xy}=K\lr{X_1^{1/p^\infty},Y_1^{1/p^\infty},X_1^{-1/p^\infty},Y_1^{-1/p^\infty} }\text{ and  }U''_{xy}=K\lr{X_2^{1/p^\infty},Y_2^{1/p^\infty},X_2^{-1/p^\infty},Y_2^{-1/p^\infty} }
\end{equation}

Then the blow up $Z$ is given as 
\begin{equation}
Z=K\lr{X_1^{1/p^\infty},Y_1^{1/p^\infty} }\bigcup_{U'_{xy}\simeq U''_{xy}} K\lr{X_2^{1/p^\infty},Y_2^{1/p^\infty} }
\end{equation}

The identification is given by the mapping $X_1\mapsto X_2Y_2$ and $Y_1\mapsto X_2^{-1}$

\subsection{Elementary Description}
Let $A=K\lr{X_0^{1/p^\infty}, \ldots,X_n^{1/p^\infty}}$ be the ring and $I$ be an ideal of the ring. We want to take fractional powers of $I$ and we do not have a general analogue on commutative algebra. Thus, we make an assumption with $I=\lr{f_1,\ldots, f_n}$ all with the same degree and say it is possible to take $p$th power roots of $f_i$ in the sense of using newton's binomial theorem. Thus, our ideal $I$ now immediately becomes something generated by infinitely many polynomials of the form $f_i^af_j^b$ where $a+b$ is same as the degree of each $f_i\in I$. This problem can be resolved if we only consider localization at $f_i$ (rather than the product $f_i^af_j^b$).

We can thus consider graded $A$ algebra with grading done by degree $d$. 

\begin{equation}
\tilde{A}=\oplus_{d\in \ds{Z}[1/p],d\geq 0}I^d \text{ where }I^0:=A
\end{equation}

Following \cite[pp 318]{liu2002algebraic}, let $t_i\in I=\tilde{A}_1$ and we have a surjective homomorphism given as 

\begin{equation}
\phi:A[T_1^{1/p^\infty},\ldots, T_n^{1/p^\infty}]\ra\tilde{A}~\text{ where }~T_i\mapsto t_i
\end{equation}
We are mostly interested in $\prj \tilde{A}$

\subsection{Universal Property of Blow UP}

We will follow \cite[p 603]{ravi2} to give a universal definition of blow up for $X\hookrightarrow\perfp$. We are considering subvariety of $\perfp$ which is locally (on affine patches) defined by a single regular element.

The Blow Up of $X\hookrightarrow\perfp$ is a cartesian diagram
\begin{equation}
\begin{tikzpicture}[descr/.style={fill=white,inner sep=1.5pt}]
        \matrix (m) [
            matrix of math nodes,
            row sep=2.5em,
            column sep=2.5em,
            text height=1.5ex, text depth=0.25ex
        ]
        { |[name=a1]|  E_X\perfp & |[name=a2]|Bl_X\perfp \\
            |[name=b1]|X & |[name=b2]|\perfp \\
                    };

         \draw [ right hook-latex]
        (a1) edge (a2)
        (b1) edge (b2)
        ;
   \draw [->]       
(a1) edge (b1)
(a2) edge node [xshift=1ex] {$\beta$}(b2)

;
\end{tikzpicture}
\end{equation}
such that $E_X\perfp$ is an effective Cartier divisor (locally generated by a single element) on $Bl_X\perfp$, where $E_X\perfp$ is pull back of $X$ by $\beta$ and the above is universal. If there is any other diagram with $D$ an effective cartier divisor on $W$ then it will uniquely factor through it

\begin{equation}
\begin{tikzpicture}[descr/.style={fill=white,inner sep=1.5pt}]
        \matrix (m) [
            matrix of math nodes,
            row sep=2.5em,
            column sep=2.5em,
            text height=1.5ex, text depth=0.25ex
        ]
        { |[name=a1]|  D & |[name=a2]|W \\
            |[name=b1]|X & |[name=b2]|\perfp \\
                    };

         \draw [ right hook-latex]
        (a1) edge (a2)
        (b1) edge (b2)
        ;
   \draw [->]       
(a1) edge (b1)
(a2) edge (b2)

;
\end{tikzpicture}
\hspace{1cm}
\begin{tikzpicture}[descr/.style={fill=white,inner sep=1.5pt}]
        \matrix (m) [
            matrix of math nodes,
            row sep=2.5em,
            column sep=2.5em,
            text height=1.5ex, text depth=0.25ex
        ]
        { |[name=c1]|  D & |[name=c2]|W \\
        |[name=a1]|  E_X\perfp & |[name=a2]|Bl_X\perfp \\
            |[name=b1]|X & |[name=b2]|\perfp \\
                    };

         \draw [ right hook-latex]
         (c1) edge (c2)
        (a1) edge (a2)
        (b1) edge (b2)
        ;
   \draw [->]  
 (c1) edge (a1) 
  (c2) edge (a2)       
(a1) edge (b1)
(a2) edge (b2)

;
\end{tikzpicture}
\end{equation}

Clearly Blowup is unique upto isomorphism and $X$ will be called the center of Blowup. If $X$ is a effective cartier divisor then, the universal property of blow up tells us that $\bl_X\perfp\ra\perfp$ is an isomorphism.

\subsection{Constructing the Blowup}

The blowup in algebraic geometry is built via $\curly{I}^d$ where $\curly{I}$ is the ideal corresponding to $X$ and $d$ represents ideal multiplication. In our case we use $d$ as a degree where $d\in\ds{Z}[1/p]$. Thus, the existence of blow up would depend upon if we can take rational powers (of the form $a/p^b$) of polynomials. The binomial theorem with rational powers can be used to construct individual examples.

In this section we want to show that the blow up corresponding to $X\ra\perfp$ is
\begin{equation}
\prj_{\ds{P}^n}(\curly{O}_{\perfp}\oplus_{d>0}I^d) \text{ where }d\in\ds{Z}[1/p]
\end{equation}
where $I$ is the ideal corresponding to $X$. We also want $I$ to be finitely generated so that we can form an compact affine covering (which is needed in the proof). We want to consider blow up as a subsheaf of the sheaf
\begin{equation}
\oplus_{d\geq0}t^d\curly{O}_{\perfp}, \text{ where $d$ is degree and }d\in\ds{Z}[1/p]
\end{equation}
in the sense given on \cite[pp 171]{eisenbud2006geometry} and the exceptional divisor is 
\begin{equation}
\prj\left(i^*\bigoplus_{d\geq 0}\curly{I}^d\right)\text{ for }i:X\ra\perfp
\end{equation}

Proof following \cite[pp 410]{gortz2010algebraic}: We work with the case of affine target with ideal $I\subset A$ (where $A$ corresponds to $\perfp$) and the principal effective Cartier divisor is $f$ on the affine source. Let $B=\oplus I^d$ and $A[If^{-1}]$ be generated by elements of the form $x/f$ where $x\in I$, such an element might be considered as degree zero element of $B_f$. Thus we have an isomorphism
\begin{equation}
A[If^{-1}]\ra B_{(f)}
\end{equation}
with inverse given as $y\in I^d$ getting mapped to $y/f^d$ (we see here that the assumption $f^d$ needs to make sense for $d$ a rational number). As $f$ runs through generating set of $I$ (we need finite generation here) the sets $D_+(f)$ form an open affine covering of the blow up. Thus, we can talk about the map
\begin{equation}
A[If^{-1}]=\frac{A[(T_\alpha)_\alpha]}{(fT_\alpha=x_\alpha)_\alpha}\ra B_{(f)}\text{ where }(x_\alpha)_\alpha\text{ is a generating set of } I
\end{equation}

which agrees with morphism above away from $V(f)$.

For example, blowing up $\ds{A}^2$ at the origin gives us \cite[p 603]{ravi2} the blow up below. Note that $x,y$ have degree $0$ and $X,Y$ have degree one.
\begin{equation}
B=K\lr{x^{1/p^\infty},y^{1/p^\infty}}, I=(x,y)\text {and } \prj(B\oplus_{d\geq 0} I^d)=\prj (B\lr{X^{1/p^\infty},Y^{1/p^\infty}})
\end{equation}

\subsection{Normal Cone} For a closed subscheme $X$ in $\perfp$ or $\arov{n}$ cut out by $I$, the magic of braided dimension helps define the normal cone as a tuple below.
\begin{equation}
\begin{aligned}
&{\curly{O}}/{I}\oplus{I}/{I^2}\oplus {I^2}/{I^3}\oplus\cdots,\\
&{\curly{O}}/{I^{1/p}}\oplus{I^{1/p}}/{I^{2/p}}\oplus {I^{2/p}}/{I^{3/p}}\oplus\cdots,\\
&{\curly{O}}/{I^{1/p^2}}\oplus{I^{1/p^2}}/{I^{2/p^2}}\oplus {I^{2/p^2}}/{I^{3/p^2}}\oplus\cdots,\\
& \qquad\qquad\qquad\vdots\qquad\\
&{\curly{O}}/{I^{1/p^i}}\oplus{I^{1/p^i}}/{I^{2/p^i}}\oplus {I^{2/p^i}}/{I^{3/p^i}}\oplus\cdots,\\
& \qquad\qquad\qquad\vdots\qquad\\
\end{aligned}
\end{equation}

As usual, the first co-ordinate is the definition in algebraic geometry.

\bibliographystyle{apalike}
\bibliography{\myreferences}

\begin{thebibliography}{}

\bibitem[Bedi, 2018]{Bedi2018}
Bedi, H. (2018).
\newblock {\em Line bundles of Rational Degree on Perfectoid Spaces}.
\newblock {PhD} dissertation, George Washington University.

\bibitem[Das, 2016]{Das_2016}
Das, S. (2016).
\newblock {\em Vector Bundles on Perfectoid Spaces}.
\newblock {PhD} dissertation, UC San Diego.

\bibitem[Eisenbud and Harris, 2006]{eisenbud2006geometry}
Eisenbud, D. and Harris, J. (2006).
\newblock {\em The Geometry of Schemes}.
\newblock Graduate Texts in Mathematics. Springer New York.

\bibitem[G{\"o}rtz and Wedhorn, 2010]{gortz2010algebraic}
G{\"o}rtz, U. and Wedhorn, T. (2010).
\newblock {\em Algebraic Geometry: Part I: Schemes. With Examples and
  Exercises}.
\newblock Advanced Lectures in Mathematics. Vieweg+Teubner Verlag.

\bibitem[Hardy and Ramanujam, 1918]{ram1}
Hardy, G. and Ramanujam, S. (1918).
\newblock Asymptotic formuae in combinatory analysis.
\newblock {\em Proceedings of the London Mathematical Society}, XVII:75--115.

\bibitem[Hartshorne, 1977]{hartshorne1977algebraic}
Hartshorne, R. (1977).
\newblock {\em Algebraic Geometry}.
\newblock Encyclopaedia of mathematical sciences. Springer.

\bibitem[Hatcher, 2002]{hatcher2002algebraic}
Hatcher, A. (2002).
\newblock {\em Algebraic Topology}.
\newblock Algebraic Topology. Cambridge University Press.

\bibitem[{Kedlaya}, 2017]{Kiran_AWS}
{Kedlaya}, K.~S. (2017).
\newblock {Sheaves, stacks, and shtukas: Arizona Winter School}.
\newblock \url{http://swc.math.arizona.edu/aws/2017/2017KedlayaNotes.pdf}.

\bibitem[{Kedlaya} and {Liu}, 2015]{Kiran_Liu_1}
{Kedlaya}, K.~S. and {Liu}, R. (2015).
\newblock {Relative p-adic Hodge theory: Foundations}.
\newblock {\em Asterisque}.

\bibitem[Koll{\'a}r, 2000]{kollar2000complex}
Koll{\'a}r, J. (2000).
\newblock {\em Complex Algebraic Geometry}.
\newblock IAS Park City mathematics series: Institute for Advanced Study.
  American Mathematical Society, Institute for Advanced Study.

\bibitem[Liu, 2002]{liu2002algebraic}
Liu, Q. (2002).
\newblock {\em Algebraic Geometry and Arithmetic Curves}.
\newblock Oxford graduate texts in mathematics. Oxford University Press.

\bibitem[Scholze, 2012]{scholze_1}
Scholze, P. (2012).
\newblock Perfectoid spaces.
\newblock {\em Publ. math. IHES}, 116:245–313.
\newblock
  \url{http://math.stanford.edu/~conrad/Perfseminar/refs/perfectoid.pdf}.

\bibitem[{Stacks Project Authors}, 2016]{stacks-project}
{Stacks Project Authors}, T. (2016).
\newblock \itshape stacks project.
\newblock \url{http://stacks.math.columbia.edu}.

\bibitem[Vakil, 2017]{ravi2}
Vakil, R. (2017).
\newblock Foundations of algebraic geometry.
\newblock \url{http://math.stanford.edu/~vakil/216blog/FOAGfeb0717public.pdf}.

\end{thebibliography}
\end{document}